\documentclass[12pt]{amsart}
\usepackage{amssymb,amsmath,amsfonts,latexsym}
\usepackage{bm}

\usepackage{array,graphics,color}
\newcommand{\newsection}[1]{\setcounter{equation}{0} \section{#1}}
\numberwithin{equation}{section}

\newtheorem{propn}{Proposition}[section]
\newtheorem{thm}[propn]{Theorem}
\newtheorem{lemma}[propn]{Lemma}

\newtheorem*{thm*}{Theorem}
\theoremstyle{definition}

\newtheorem{rem}{Remark}[section]

\newcommand{\Hil}{\mathcal{H}}

\newcommand{\Nat}{\mathbb{N}}

 \newcommand{\D}{\mathbb{D}}

\newcommand{\clb}{\mathcal{B}}

\newcommand{\cld}{\mathcal{D}}
\newcommand{\cle}{\mathcal{E}}

\newcommand{\clh}{\mathcal{H}}
\newcommand{\clk}{\mathcal{K}}

\newcommand{\clq}{\mathcal{Q}}

\newcommand{\cls}{\mathcal{S}}

\newcommand{\raro}{\rightarrow}
\newcommand{\NI}{\noindent}

\begin{document}

\title{Ando dilations, von Neumann inequality, and distinguished varieties}

\author[Das] {B. Krishna Das}
\address{Department of Mathematics, Indian Institute of Technology Bombay, Powai, Mumbai, 400076, India}
\email{dasb@math.iitb.ac.in, bata436@gmail.com}

\author[Sarkar]{Jaydeb Sarkar}
\address{Indian Statistical Institute, Statistics and Mathematics Unit, 8th Mile, Mysore Road, Bangalore, 560059, India}
\email{jay@isibang.ac.in, jaydeb@gmail.com}


\begin{abstract}
Let $\mathbb{D}$ denote the unit disc in the complex plane
$\mathbb{C}$ and let $\mathbb{D}^2 = \mathbb{D} \times \mathbb{D}$
be the unit bidisc in $\mathbb{C}^2$. Let $(T_1, T_2)$ be a pair of
commuting contractions on a Hilbert space $\clh$. Let $\dim
{\textup{ran}}(I_{\mathcal{H}} - T_j T_j^*) < \infty$, $j = 1, 2$,
and let $T_1$ be a pure contraction. Then there exists a variety $V
\subseteq \overline{\D}^2$ such that for any polynomial $p \in
\mathbb{C}[z_1, z_2]$, the inequality
\[
\|p(T_1,T_2)\|_{\mathcal{B}(\mathcal{H})} \leq \|p\|_V
\]
holds. If, in addition, $T_2$ is pure, then
\[V = \{(z_1, z_2) \in
\mathbb{D}^2: \det (\Psi(z_1) - z_2 I_{\mathbb{C}^n}) = 0\}\]is a
distinguished variety, where $\Psi$ is a matrix-valued analytic
function on $\mathbb{D}$ that is unitary on $\partial \mathbb{D}$.
Our results comprise a new proof, as well as a generalization, of
Agler and McCarthy's sharper von Neumann inequality for pairs of
commuting and strictly contractive matrices.

\end{abstract}

\maketitle

\section*{Notation}

\begin{list}{\quad}{}

\item $\mathbb{D}$ \quad \quad \quad\quad Open unit disc in the complex plane
$\mathbb{C}$.

\item  $\mathbb{D}^2$ \; \quad \quad\quad Open unit bidisc in $\mathbb{C}^2$.

\item  $\clh$, $\cle$ \; \quad\quad Hilbert spaces.

\item $\clb(\clh)$ \quad \;\quad The space of all bounded linear
operators on $\clh$.

\item $H^2_{\cle}(\mathbb{D})$ \quad \quad  $\cle$-valued Hardy space on
$\mathbb{D}$.

\item $M_z$  \quad \, \quad \quad Multiplication operator by the coordinate function
$z$.

\item $H^\infty_{\clb(\cle)}(\D)$ \quad Set of $\clb(\cle)$-valued bounded analytic functions on
$\mathbb{D}$.

\end{list}

All Hilbert spaces are assumed to be over the complex numbers. For a
closed subspace $\cls$ of a Hilbert space $\clh$, we denote by
$P_{\cls}$ the orthogonal projection of $\clh$ onto $\cls$.

\newsection{Introduction}

The famous von Neumann inequality \cite{vN} states that: if $T$ is a
linear operator on a Hilbert space $\clh$ of norm one or less (that
is, $T$ is a contraction), then for any polynomial $p \in
\mathbb{C}[z]$, the inequality
\[\|p(T)\|_{\clb(\clh)} \leq \|p\|_{\mathbb{D}}\]holds. Here $\|p\|_{\mathbb{D}}$
denotes the supremum of $|p(z)|$ over the unit disc $\D$.

In 1953, Sz.-Nagy \cite{Nagy} proved that a linear operator on a
Hilbert space is a contraction if and only if the operator has a
unitary dilation. This immediately gives a simple and elegant proof
of the von Neumann inequality.

\vspace{0.1in}

In 1963, Ando \cite{Ando} proved the following generalization of
Sz.-Nagy's dilation theorem: Any pair of commuting contractions has
a commuting unitary dilation. As an immediate consequence, we obtain
the following two variables von Neumann inequality:

\NI\textsf{Theorem (Ando):} Let $(T_1, T_2)$ be a pair of commuting
contractions on a Hilbert space $\clh$. Then for any polynomial $p
\in \mathbb{C}[z_1, z_2]$, the inequality
\[\|p(T_1, T_2)\|_{\clb(\clh)} \leq \|p\|_{\D^2}\]holds.

However, for three or more commuting contractions the above von
Neumann type inequality is not true in general (see \cite{CD},
\cite{V}). An excellent source of further information on von Neumann
inequality is the monograph by Pisier \cite{Pi}.

In a recent seminal paper, Agler and McCarthy \cite{AM1} proved a
sharper version of von-Neumann inequality for pairs of commuting and
strictly contractive matrices (see Theorem 3.1 in \cite{AM1}): Two
variables von Neumann inequality can be improved in the case of a
pair of commuting and strictly contractive operators $(T_1, T_2)$ on
a finite dimensional Hilbert space $\clh$ to
\[\|p(T_1, T_2)\|_{\clb(\clh)} \leq \|p\|_{V} \quad \quad (p \in \mathbb{C}[z_1, z_2]),\]
where V is a distinguished variety depending on the pair $(T_1,
T_2)$. Again, $\|p\|_V$ is the supremum of $|p(z_1, z_2)|$ over $V$.
The proof of this result involves many different techniques
including isometric dilation of $(T_1, T_2)$ to a vector-valued
Hardy space and approximation of commuting matrices by digonalizable
commuting matrices.

We recall that a non-empty set $V$ in $\mathbb{C}^2$ is a
\textit{distinguished variety} if there is a polynomial $p \in
\mathbb{C}[z_1, z_2]$ such that
\[V = \{(z_1, z_2) \in \mathbb{D}^2 : p(z_1, z_2) = 0\},\]and
$V$ exits the bidisc through the distinguished boundary, that
is, \[\overline{V} \cap \partial \mathbb{D}^2 = \overline{V} \cap
(\partial \mathbb{D} \times \partial \mathbb{D}).\]Here $\partial
\mathbb{D}^2$ and $\partial\mathbb{D}\times \partial\mathbb{D}$ denote the
boundary and the distinguished boundary of the bidisc repectively, and
$\overline{V}$ is the closure of $V$ in $\overline{\mathbb{D}^2}$.
We denote by $\partial V$ the set $\overline{V} \cap \partial
\mathbb{D}^2$, the boundary of $V$ within the zero set of the
polynomial $p$ and $\overline{\mathbb{D}^2}$.

In the same paper \cite{AM1}, Agler and McCarthy proved that a
distinguished variety can be represented by a rational matrix inner
function in the following sense (see Theorem 1.12 in \cite{AM1}):
Let $V \subseteq \mathbb{C}^2$. Then $V$ is a distinguished variety
if and only if there exists a rational matrix inner function $\Psi
\in H^\infty_{\clb(\mathbb{C}^n)}(\D)$, for some $n \geq 1$, such
that
\[V = \{(z_1, z_2) \in \mathbb{D}^2: \det (\Psi(z_1) - z_2 I_{\mathbb{C}^n}) =
0\}.\]The proof uses dilation and model theoretic techniques (see
page 140 in \cite{AM1}) in the sense of Sz.-Nagy and Foias
\cite{NF}. See Knese \cite{K} for another proof.

In this paper, in Theorem \ref{pure isometric dilation}, we obtain
an explicit way to construct isometric dilations of a large class of
commuting pairs of contractions: Let $(T_1, T_2)$ be a pair of
commuting operators on a Hilbert space $\clh$ and $\|T_j\| \leq 1$,
$j = 1, 2$. Let $\dim {\textup{ran}}(I_{\clh} - T_j T_j^*) <
\infty$, $j = 1, 2$, and $T_1$ be a pure contraction (that is,
$\lim_{m \raro \infty} \|T_1^{*m} h\| =0$ for all $h \in \clh$). Set
$\cle = \mbox{ran}(I_{\clh} - T_1 T_1^*)$. Then there is an
$\mathcal{B}(\cle)$-valued inner function $\Psi$ such that the commuting isometric pair
$(M_z, M_{\Psi})$ on the $\cle$-valued Hardy space
$H^2_{\cle}(\mathbb{D})$ is an isometric dilation of $(T_1, T_2)$.

\NI We actually prove a more general dilation result in Theorem
\ref{gen-dilation}.

Then in Theorem \ref{vn-new1}, we prove: There exists a variety $V
\subseteq \overline{\D}^2$ such that
\[
 \|p(T_1,T_2)\|_{\clb(\clh)} \le \|p\|_V \quad
\quad (p \in \mathbb{C}[z_1, z_2]).
\]
If, in addition, $T_2$ is pure, then $V$ can be taken to be a
distinguished variety.

It is important to note that every distinguished variety, by
definition, is a subset of the bidisc $\D^2$.

Our results comprise both a new proof, as well as a generalization,
of Agler and McCarthy's sharper von Neumann inequality for pairs of
commuting and strictly contractive matrices  (see the final
paragraph in Section \ref{S-vne}).


The remainder of this paper is built as follows. In Section
\ref{prel}, we first recall some basic definitions and results. We
then proceed to prove an important lemma which will be used in the
sequel. Dilations of pairs of commuting contractions are studied in
Section \ref{S-d}. In Section \ref{S-vne}, we use results from the
previous section to show a sharper von Neumann inequality for pairs
of pure commuting contractive tuples with finite dimensional defect spaces.
In the concluding section,
Section \ref{S-c}, among other things, we prove that the
distinguished variety in our von Neumann inequality is independent
of the choice of $(T_1, T_2)$ and $(T_2, T_1)$.

\newsection{Preliminaries and a Correlation lemma}\label{prel}

First we recall some definitions of objects we are going to use and
fix few notations.

Let $T$ be a contraction on a Hilbert space $\clh$ (that is, $\|T f
\| \leq \|f\|$ for all $f \in \clh$ or, equivalently, if $I_{\clh} -
T T^* \geq 0$). Recall again that $T$ is \textit{pure} if $\lim_{m
\raro 0} \|T^{*m} f\| =0$ for all $f \in \clh$.

Let $T$ be a contraction and $\cle$ be a Hilbert space. An isometry
$\Gamma : \clh \raro H^2_{\cle}(\D)$ is called an isometric dilation
of $T$ (cf. \cite{JS}) if \[\Gamma T^* = M_z^* \Gamma.\] If, in
addition,
\[H^2_{\cle}(\D) = \overline{\mbox{span}} \{z^m \Gamma f : m \in
\mathbb{N}, f \in \clh\},\]then we say that $\Gamma : \clh \raro
H^2_{\cle}(\D)$ is a minimal isometric dilation of $T$.

Now let $T$ be a pure contraction on a Hilbert space $\clh$. Set
\[\cld_T = \overline{\mbox{ran}}(I_{\clh} - T T^*), \quad \quad D_T =
(I_{\clh} - T T^*)^{\frac{1}{2}}.\]
Then $\Pi : \clh \raro H^2_{\cld_T}(\D)$ is a minimal isometric
dilation of $T$ (cf. \cite{JS}), where
\begin{equation}\label{dil-def} (\Pi h)(z) = D_T(I_{\clh} - z
T^*)^{-1}h \quad \quad (z \in \D, h \in \clh).\end{equation} In
particular, $\clq := \mbox{ran} \Pi$ is a $M_z^*$-invariant subspace
of $H^2_{\cld_T}(\D)$ and \[T \cong P_{\clq} M_z|_{\clq}.\]

Finally, recall that a contraction $T$ on $\clh$ is said to be
\textit{completely non-unitary} if there is no non-zero $T$-reducing
subspace $\cls$ of $\clh$ such that $T|_{\cls}$ is a unitary
operator. It is well known that for every contraction $T$ on a
Hilbert space $\clh$ there exists a unique \textit{canonical
decomposition} $\clh = \clh_0 \oplus \clh_1$ of $\clh$ reducing $T$,
such that $T|_{\clh_0}$ is unitary and $T|_{\clh_1}$ is completely
non-unitary. We therefore have the following decomposition of $T$:
\[T = \begin{bmatrix} T|_{\clh_0} & 0 \\ 0 &
T|_{\clh_1}\end{bmatrix}.\]

We now turn to the study of contractions with finite dimensional
defect spaces. Let $(T_1, T_2)$ be a pair of commuting contractions
and $\dim \cld_{T_j} < \infty$, $j = 1, 2$. Since \[(I_{\clh} - T_1
T_1^*) + T_1 (I_{\clh} - T_2 T_2^*) T_1^* = T_2 (I_{\clh} - T_1
T_1^*) T_2^* + (I_{\clh} - T_2 T_2^*),\] it follows that \[
 \|D_{T_1}h\|^2+ \|D_{T_2}T_1^*h\|^2=\|D_{T_1}T_2^*h\|^2+\|D_{T_2}h\|^2
 \quad (h\in\Hil).
\]Thus \[U : \{D_{T_1} h \oplus D_{T_2} T_1^* h : h \in \clh\} \raro
\{D_{T_1} T_2^* h \oplus D_{T_2} h : h \in \clh\}\]defined by
\begin{equation}\label{U-h}U\left(D_{T_1}h, D_{T_2}T_1^*
h\right)=\left(D_{T_1} T_2^*h, D_{T_2}h\right) \quad \quad (h \in
\clh),\end{equation}is an isometry. Moreover, since $\mbox{dim~}
\cld_{T_j} < \infty$, $j = 1, 2$, it follows that $U$ extends to a
unitary, denoted again by $U$, on $\cld_{T_1} \oplus \cld_{T_2}$. In
particular, there exists a unitary operator
\begin{equation}
\label{unitary} U: = \begin{bmatrix}A&B\\C&D\end{bmatrix}:
\cld_{T_1}\oplus\cld_{T_2}\to \cld_{T_1}\oplus \cld_{T_2},
\end{equation}
such that (\ref{U-h}) holds.

The following lemma plays a key role in our considerations.

\begin{lemma}\label{identity}
Let $(T_1, T_2)$ be a pair of commuting contractions on a Hilbert
space $\clh$. Let $T_1$ be pure and $\mbox{dim~} \cld_{T_j} <
\infty$, $j = 1, 2$. Then with $U =
\begin{bmatrix}A&B\\C&D\end{bmatrix}$ as above we have
\[
D_{T_1}T_2^* = AD_{T_1} + \sum_{n= 0}^\infty B D^{n}CD_{T_1} T_1^{*
n+1},
\]
where the series converges in the strong operator topology.
\end{lemma}
\NI\textsf{Proof.} For each $h \in \clh$ we have
\[\begin{bmatrix}A&B\\C&D\end{bmatrix} \begin{bmatrix} D_{T_1}h\\
D_{T_2}T_1^*h \end{bmatrix}= \begin{bmatrix} D_{T_1}T_2^{*}h\\ D_{T_2}h\end{bmatrix},\]
that is,
\begin{equation} \label{one} D_{T_1} T_2^{*} h= AD_{T_1}h+
BD_{T_2}T_1^*h,
\end{equation}
and
\begin{equation}\label{two}
D_{T_2}h=CD_{T_1}h+ DD_{T_2}T_1^* h.
\end{equation}
By replacing $h$ by $T_1^{*}h$ in ~\eqref{two}, we have
\[D_{T_2}T_1^*h = C D_{T_1} T_1^* h+ DD_{T_2}T_1^{* 2} h.\]
Then by ~\eqref{one}, we have
\[\begin{split}
D_{T_1}T_2^*h & =AD_{T_1}h + BD_{T_2}T_1^*h  \\ & = AD_{T_1}h + B (C
D_{T_1} T_1^* h+ DD_{T_2}T_1^{* 2} h), \end{split}\] that is,
\begin{equation}\label{three}D_{T_1}T_2^*h = AD_{T_1}h+ B C
D_{T_1}T_1^*h + BDD_{T_2}T_1^{* 2}h.\end{equation} In \eqref{three},
again replacing $D_{T_2}T_1^{* 2}h$ by $CD_{T_1}T_1^{* 2}h+
DD_{T_2}T_1^{* 3}h$, we have
\[D_{T_1}T_2^*h = AD_{T_1}h+ B C
D_{T_1}T_1^*h + BD C D_{T_1}T_1^{* 2}h + B D^2 D_{T_2} T_1^{* 3}h.\]
Going on in this way, we obtain
\[
D_{T_1}T_2^*h = AD_{T_1}h+  \sum_{n= 0}^m B D^{n}CD_{T_1} T_1^{*
n+1} h +  (B D^{m+1} D_{T_2} T_1^{* m+2} h)\quad (m\in\Nat).
\]
By $\|D\| \leq 1$ and \[\lim_{m \raro \infty} T_1^{*m} h = 0,\] it
follows that \[\lim_{m \raro \infty} (B D^m D_{T_2} T_1^{* m+1} h) =
0.\] Finally, the convergence of the desired series follows from the
fact that \[\begin{split}\|D_{T_1} T_2^* h - A D_{T_1} h -
\sum_{n=0}^m B D^{n}CD_{T_1} T_1^{* n+1} h\| & = \|B D^{m+1}CD_{T_2}
T_1^{* m+2} h\| \\ & \leq \|T_1^{* m+2} h\|,\end{split}\]for all $m
\in \mathbb{N}$ and again that $\lim_{m \raro \infty} T_1^{*m} h =
0$. This completes the proof of the lemma. \qed

\begin{rem}\label{inf-D}
With the above assumptions, the conclusion of Lemma \ref{identity}
remains valid if we add the possibility \[\dim \cld_{T_1} = \infty,
\;\; or \;\; \dim \cld_{T_2} = \infty.\] To this end, let $\dim
\cld_{T_1} = \infty$, or $\dim \cld_{T_2} = \infty$. Then there
exists an infinite dimensional Hilbert space $\cld$ such that the
isometry
\[U : \{D_{T_1} h \oplus D_{T_2} T_1^* h : h \in
\clh\} \oplus \{0_{\cld}\} \raro \{D_{T_1} T_2^* h \oplus D_{T_2} h
: h \in \clh\} \oplus \{0_{\cld}\}\]defined by
\[U\left(D_{T_1}h, D_{T_2}T_1^* h , 0_{\cld} \right)=\left(D_{T_1}
T_2^*h, D_{T_2}h , 0_{\cld}\right) \quad
\quad (h \in \clh),\] extends to a unitary, denoted again by $U$, on
$\cld_{T_1} \oplus (\cld_{T_2} \oplus \cld)$. Now we proceed
similarly with the unitary matrix \[U =
\begin{bmatrix}A&B\\C&D\end{bmatrix} \in \clb(\cld_{T_1} \oplus
(\cld_{T_2} \oplus \cld)),\] to obtain the same conclusion as in
Lemma \ref{identity}.
\end{rem}

\newsection{Ando Dilations}\label{S-d}

One of the aims of this section is to obtain an explicit isometric
dilation of a pair of commuting tuple of pure contractions (see
Theorem \ref{gen-dilation}). In particular, the construction of
isometric dilations of pairs of commuting contractions with finite
defect indices (see Theorem \ref{pure isometric dilation}) will be
important to us in the sequel.

We begin by briefly recalling some standard facts about transfer
functions (cf. \cite{AM-b}). Let $\clh_1$ and $\clh_2$ be two
 Hilbert spaces, and
\[U = \begin{bmatrix}A&B\\C&D\end{bmatrix} \in \clb(\clh_1 \oplus
\clh_2),\]be a unitary operator. Then the $\clb(\clh_1)$-valued
analytic function $\tau_U$ on $\mathbb{D}$ defined by
\[\tau_U (z) := A + z B (I-z D)^{-1} C \quad \quad (z \in \D), \]
is called the \textit{transfer function} of $U$. Using $U^* U = I$,
a standard computation yields (cf. \cite{AM-b})
\begin{equation}
\label{isometry}
 I- \tau_U (z)^*\tau_U (z)=(1-|z|^2) C^*(I-\bar{z} D^*)^{-1}(I-z D)^{-1} C \quad \quad (z\in\D).
\end{equation}

Now let $\clh_1 = \mathbb{C}^m$ and $\clh_2 = \mathbb{C}^n$, and let
$U$ be as above. Then $\tau_U$ is a contractive rational
matrix-valued function on $\D$. Moreover, $\tau_U$ is unitary on
$\partial \mathbb{D}$ (see page 138 in \cite{AM1}). Thus $\tau_U$ is
a rational matrix-valued inner function.

We now turn our attention to the study of the transfer function of
the unitary matrix $U^*$ in ~\eqref{unitary}. Set
\begin{equation}\label{psi}\Psi(z) : = \tau_{U^*}(z) = A^*+ zC^*(I-zD^*)^{-1}B^* \quad \quad
(z\in\D).\end{equation} Then $\Psi$ is a
$\mathcal{B}(\cld_{T_1})$-valued inner function on $\D$. Thus the
multiplication operator $M_{\Psi}$ on $H^2_{\cld_{T_1}}(\D)$ defined
by
\[(M_{\Psi} f)(w) = \Psi(w) f(w) \quad \quad (w \in \D,
f \in H^2_{\cld_{T_1}}(\D)),\]
is an isometry (cf. \cite{NF}).

Now, we are ready to prove our first main result.

\begin{thm}\label{pure isometric dilation}
Let $(T_1, T_2)$ be a pair of commuting contractions on a Hilbert
space $\clh$. Let $T_1$ be pure and $\mbox{dim~} \cld_{T_j} <
\infty$, $j = 1,2$. Then there exists an isometry $\Pi : \clh \raro
H^2_{\cld_{T_1}}(\D)$ and an inner function $\Psi \in
H^\infty_{\clb(\cld_{T_1})}(\D)$ such that
\[\Pi T_1^* = M_z^* \Pi,\] and
\[\Pi T_2^* = M_{\Psi}^* \Pi.\] Moreover
\[T_1\cong P_{\clq} M_z|_{\clq}, \quad \mbox{and} \quad\quad
T_2\cong P_{\clq}M_{\Psi}|_{\clq},\]where $\clq := \mbox{ran~} \Pi$
is a joint $(M_z^*, M_{\Psi}^*)$-invariant subspace of
$H^2_{\cld_{T_1}}(\D)$.
\end{thm}
\NI\textsf{Proof.} Let $\Pi : \clh \raro H^2_{\cld_{T_1}}(\D)$ be
the minimal isometric dilation of $T_1$ as defined in
(\ref{dil-def}) and let $\Psi$ be as in (\ref{psi}). Now it is
enough to show that $\Pi$ intertwine $T_2^*$ and $M_{\Psi}^*$. Let
$h \in \clh$, $n \geq 0$ and $\eta \in \cld_{T_1}$. Then
\[\begin{split}
\langle M_{\Psi}^* \Pi h, z^n \eta \rangle & = \langle \Pi h, \Psi
(z^n \eta) \rangle \\& = \langle D_{T_1} (I_{\clh} - z T_1^*)^{-1}
h, (A^* + C^* \sum_{q=0}^\infty D^{*q} B^* z^{q+1}) (z^n \eta)
\rangle
\\& = \langle D_{T_1} \sum_{p=0}^\infty z^p T_1^{*p} h, (A^* + C^*
\sum_{q=0}^\infty D^{*q} B^* z^{q+1}) (z^n \eta) \rangle \\& =
\langle D_{T_1} T_1^{*n} h, A^* \eta \rangle + \sum_{q=0}^\infty
\langle D_{T_1} T_1^{* q+n+1}h, C^* D^{*q} B^* \eta \rangle \\& =
\langle A D_{T_1} T_1^{*n} h, \eta \rangle + \sum_{q=0}^\infty
\langle B D^q C D_{T_1} T_1^{* q+n+1}h, \eta \rangle \\& = \langle
(A D_{T_1} + \sum_{q=0}^\infty B D^q C D_{T_1} T_1^{* q +1})
T_1^{*n} h, \eta \rangle.
\end{split}\]
Then by Lemma \ref{identity}, we have
\[\langle M_{\Psi}^* \Pi h, z^n
\eta \rangle = \langle D_{T_1} T_2^* (T_1^{*n}h), \eta \rangle.
\]
Hence
\[\begin{split}
\langle \Pi T_2^* h, z^n \eta \rangle & = \langle D_{T_1} (I_{\clh}
- z T_1^*)^{-1} T_2^* h, z^n \eta \rangle \\ & = \langle
D_{T_1} \sum_{p=0}^\infty z^p T_1^{*p} T_2^* h, z^n \eta \rangle \\
& = \langle D_{T_1} T_1^{*n} T_2^* h, \eta \rangle  \\
& = \langle D_{T_1} T_2^* (T_1^{*n} h), \eta \rangle \\ & = \langle
M_{\Psi}^* \Pi h, z^n \eta \rangle.\end{split}\]This implies
\[M_{\Psi}^* \Pi = \Pi T_2^*.\]
The second claim is an immediate consequence of the first part. This completes the proof of the theorem.
\qed

Theorem \ref{pure isometric dilation} remains valid if we drop the
assumption that $\dim \cld_{T_j} < \infty$, $j = 1, 2$. Indeed, the
only change needed in the proof of Theorem \ref{pure isometric
dilation} is to replace the transfer function $\Psi$ in (\ref{psi})
by the transfer function of $U^*$ in Remark \ref{inf-D}. In this
case, however, the new transfer function will be a contractive
multiplier. Thus we have proved the following dilation result.

\begin{thm}\label{gen-dilation}
Let $(T_1, T_2)$ be a pair of commuting contractions on a Hilbert
space $\clh$ and let $T_1$ be a pure contraction. Then there exists
an isometry $\Pi : \clh \raro H^2_{\cld_{T_1}}(\D)$ and a
contractive multiplier $\Psi \in H^\infty_{\clb(\cld_{T_1})}(\D)$
such that
\[\Pi T_1^* = M_z^* \Pi,\] and
\[\Pi T_2^* = M_{\Psi}^* \Pi.\] In particular,
\[T_1\cong P_{\clq} M_z|_{\clq}, \quad \mbox{and} \quad\quad
T_2\cong P_{\clq}M_{\Psi}|_{\clq},\]where $\clq := \mbox{ran~} \Pi$
is a joint $(M_z^*, M_{\Psi}^*)$-invariant subspace of
$H^2_{\cld_{T_1}}(\D)$.
\end{thm}

\newsection{von Neumann inequality}\label{S-vne}

This section is devoted mostly to the study of von Neumann
inequality for the class of pure and commuting contractive pair of
operators with finite defect spaces. However, we treat this matter
in a slightly general setting. It is convenient to be aware of the results and constructions of Section \ref{S-d}.

We begin by noting the following proposition.

\begin{propn}
\label{eigenvalues} Let $U = \begin{bmatrix} A& B\\C &
D\end{bmatrix}$ be a unitary matrix on $\clh \oplus \clk$ and let
$A$ be a completely non-unitary contraction. Then for all $z \in
\D$, $\tau_U(z)$ does not have any unimodular eigenvalues.
\end{propn}
\NI\textsf{Proof.} Let $z \in \D$ and $A$ be a completely
non-unitary contraction. Suppose, by contradiction, \[(\tau_U(z)) v
= \lambda v,\]for some non-zero vector $v \in \clh$ and for some
$\lambda\in\partial\D$. Since $\tau_U(z)$ is a contraction, we have
\[(\tau_U(z))^* v =\bar{\lambda} v,\]and hence, by ~\eqref{isometry} \[C v = 0.\]
This and the definition of $\tau_U$ implies
\[
A v = (\tau_U(z)) v = \lambda v.
\]
Then $A$ has a non-trivial unitary part. This contradiction
establishes the proposition. \qed

Now we examine the role of the unitary part of the contraction $A$
to the transfer function $\tau_U$.

\begin{propn}\label{cd-U}
Let $U = \begin{bmatrix} A& B\\C & D\end{bmatrix}$ be a unitary
matrix on $\clh \oplus \clk$ and let $A =
\begin{bmatrix} W&0\\0&A'\end{bmatrix} \in \clb(\clh_{0} \oplus \clh_{1})$ be
the canonical decomposition of $A$ into the unitary part $W$ on
$\clh_0$ and the completely non-unitary part
$A'$ on $\clh_1$. Then $U'= \begin{bmatrix} A'& B \\
C|_{\clh_1} & D
\end{bmatrix}$ is a unitary operator on $\clh_1 \oplus \clk$ and \[\tau_U(z) =
\begin{bmatrix}
W & 0 \\ 0 & \tau_{U'}(z) \end{bmatrix} \in \clb(\clh_0 \oplus
\clh_1) \quad \quad (z \in \mathbb{D}).\]
\end{propn}
\NI\textsf{Proof.} Let us observe first that if $U$ is a unitary
operator, then \[A^* A + C^* C = I_{\clh},\]and \[A A^* + B B^* =
I_{\clh}.\]On account of $A^* A|_{\clh_0} = A A^*|_{\clh_0} =
I_{\clh_0}$, the first equality implies \[C^* C|_{\clh_0} = 0,\]and
hence \[\clh_0 \subseteq \ker C,\]while the second equality yields
\[\begin{split} \overline{ran} B & = \overline{ran} (B B^*) \\ & =
\overline{ran} (I_{\clh} - A A^*) \\ & = \overline{ran} (I_{\clh_1}
- A A^*|_{\clh_1}) \\ & \subseteq \clh_1.\end{split}\]
Consequently, it follows that  $U'= \begin{bmatrix} A'& B \\
C|_{\clh_1} & D
\end{bmatrix}$ is a unitary operator on $\clh_1 \oplus \clk$, and hence
\[\tau_U(z) = W \oplus \tau_{U'}(z) \quad \quad (z \in \D),\]
follows from the definition of transfer functions. \qed

We now return to the study of rational inner functions. Let $U =
\begin{bmatrix}A& B \\ C & D\end{bmatrix}$ be a unitary on $\mathbb{C}^m\oplus
\mathbb{C}^n$. Let $A = \begin{bmatrix}W& 0\\0& E \end{bmatrix}$ on
$\mathbb{C}^m = H_0\oplus H_1$ be the canonical decomposition of $A$
into the unitary part $W$ on $H_0$ and the completely non-unitary
part $E$ on $H_1$. Then by the previous proposition, we have
\[\tau_U(z) = \begin{bmatrix} \Psi_0(z) & 0 \\ 0 & \Psi_1(z) \end{bmatrix}
\in \clb(H_0 \oplus H_1) \quad \quad (z \in \D),\] where \[\Psi_0(z)
= W \quad \quad (z \in \D)\]is a $\clb(H_0)$-valued unitary
constant, for some unitary $W \in \clb(H_0)$, and \[\Psi_1(z) = E +
z B(I - z D)^{-1} C \quad \quad (z \in \D)\]is a $\clb(H_1)$-valued
rational inner function. It should be noted, however, that the
distinguished varieties corresponding to the rational inner
functions $\tau_U$ and $\Psi_1$ are the same, that is,
\[
\{(z_1,z_2)\in\D^2:\det(\tau_U(z_1)-z_2
I_{\mathbb{C}^m})=0\}=\{(z_1,z_2)\in\D^2:\det(\Psi_1(z_1)-z_2
I_{H_1})=0\}.
\]
This follows from the observation that the unitary summand $\Psi_0$
would add sheets to the variety corresponding to $\Psi_1$ of the
type $\mathbb{C} \times \{\lambda\}$, for some $\lambda \in
\partial \D$, which is disjoint from $\D^2$.
However, we want to stress here that, by
Proposition~\ref{eigenvalues} and the fact that $\Psi_1$ is unitary
on $\partial\D$,
\[\overline{V}_{\Psi_1} = \{(z_1,z_2)\in\overline{\D^2}: \det (\Psi_1(z_1) - z_2I_{H_1}) =
0\},\]where $V_{\Psi_1}$ is the distinguished variety corresponding
to the inner multiplier $\Psi_1$ and $\overline{V}_{\Psi_1}$ is the
closure of $V_{\Psi_1}$ in $\overline{\D^2}$.

We now have all the ingredients in place to prove a von Neumann type
inequality for pairs of commuting contractions.

\begin{thm}\label{vn-new1}
Let $(T_1, T_2)$ be a pair of commuting contractions on a Hilbert
space $\clh$. Let $T_1$ be pure and $\mbox{dim~} \cld_{T_j} <
\infty$, $j = 1,2$. Then there exists a variety $V \subseteq
\overline{\D^2}$ such that
\[
 \|p(T_1,T_2)\|\le \|p\|_V \quad \quad (p \in
\mathbb{C}[z_1, z_2]).
\]
If, in addition, $T_2$ is pure, then $V$ can be taken to be a
distinguished variety.
\end{thm}
\NI \textsf{Proof.} By Theorem~\ref{pure isometric dilation}, there
is a rational inner function $\Psi \in
H^\infty_{\clb(\cld_{T_1})}(\D)$ and a joint $(M_z^*,
M_{\Psi}^*)$-invariant subspace $\clq$ of $H^2_{\cld_{T_1}}(\D)$
such that \[T_1\cong P_{\clq}M_z|_{\clq}, \quad \mbox{and} \quad
T_2\cong P_{\clq}M_{\Psi}|_{\clq}.\] Here \[\Psi(z) = \tau_{U^*}(z)
= A^* + z C^* (I - z D^*)^{-1} B^*, \quad \quad (z \in \D)\]is the
transfer function of the unitary $U^* =
\begin{bmatrix}A^*& C^* \\ B^* & D^* \end{bmatrix}$ in
$\mathcal{B}(\cld_{T_1}\oplus\cld_{T_2})$ as defined in
~\eqref{unitary}. Let $A^*= \begin{bmatrix}W&0\\ 0& E^*\end{bmatrix}
\in \mathcal{B}(H_0 \oplus H_1)$ on $\cld_{T_1}= H_0 \oplus H_1$ be
the canonical decomposition of $A^*$ in to the unitary part $W$ on
$H_0$ and the completely non-unitary part $E^*$ on $H_1$. Now from
Proposition \ref{cd-U} (or the discussion following Proposition
\ref{cd-U}) it follows that
\[\Psi(z) = \begin{bmatrix} \Psi_0(z) & 0 \\ 0 & \Psi_1(z) \end{bmatrix}
\in \clb(H_0 \oplus H_1) \quad \quad (z \in \D),\] where \[\Psi_0(z)
= W \in \clb(H_0) \quad (z \in \D),\]and
\[\Psi_1(z) = E^* + z C^*(I - z
D^*)^{-1} B^* \in \clb(H_1) \quad \quad (z \in \D).\] Let us set
\[V =V_0\cup V_1,\]where
\[V_0 =\{(z_1,z_2)\in\D\times\overline{\D}: \textup{det}(\Psi_0(z_1)- z_2 I_{H_0})=0\},\] and\[
V_1 =\{(z_1,z_2)\in\D^2: \textup{det}(\Psi_1(z_1)-z_2I_{H_1})=0\}.\]
Clearly \[V_0 =\{(z_1,z_2)\in\D\times\overline{\D}: \textup{det}(W - z_2
I_{H_0})=0\} = \mathop{\cup}_{j=1}^l \D \times \{\lambda_j\}
\subseteq \D \times \partial \D,\]where $\{\lambda_j\}_{j=1}^l =
\sigma(W) \subseteq \partial \D$. Now for each $p \in
\mathbb{C}[z_1, z_2]$, we have
\begin{align*}
\|p(T_1,T_2)\|_{\clh} &=\|P_{\clq}p(M_z,M_{\Psi})|_{\clq}\|_{\clq}\\
&\le \|p(M_z, M_{\Psi})\|_{H^2_{\cld_{T_1}}(\D)}\\
&\le \|p(M_{e^{i \theta}}, M_{\Psi(e^{i
\theta})})\|_{L^2_{\cld_{T_1}}(\mathbb{T})}\\ &= \|M_{p( e^{i
\theta} I_{\cld_{T_1}}, \Psi(e^{i
\theta}))}\|_{L^2_{\cld_{T_1}}(\mathbb{T})}\\ &= \sup_{\theta} \| p(
e^{i \theta} I_{\cld_{T_1}}, \Psi(e^{i
\theta}))\|_{\clb(\cld_{T_1})} \\ & = \sup_{\theta}\lVert
p(e^{i\theta}I_{H_0},\Psi_0(e^{i\theta}))\oplus
p(e^{i\theta}I_{H_1},\Psi_1(e^{i\theta}))\rVert_{\clb(H_0) \oplus
\clb(H_1)} \\ & \leq \max \{ \sup_{\theta}\lVert
p(e^{i\theta}I_{H_0},\Psi_0(e^{i\theta}))\rVert_{\clb(H_0)}, \lVert
p(e^{i\theta}I_{H_1},\Psi_1(e^{i\theta}))\rVert_{\clb(H_1)}\}.
\end{align*}
But now, since $\Psi(e^{i\theta})$ is unitary on $ \partial \D$, for
each fixed $e^{i \theta} \in
\partial \D $ and $j = 1, 2$, we have
\[\begin{split}
\| p( e^{i \theta} I_{H_j}, \Psi_j(e^{i \theta}))\|_{\clb(H_j)} & =
\sup \{ |p(e^{i \theta}, \lambda)|: \lambda \in \sigma(\Psi_j(e^{i
\theta}))\} \\ & = \sup \{|p(e^{i \theta}, \lambda)|: \det
(\Psi_j(e^{i \theta}) - \lambda I_{H_j}) = 0\}\\ & \leq
\|p\|_{\partial V_j},
\end{split}
\]
and hence, by continuity, we obtain \[\|p(T_1,T_2)\|_{\clh} \le \|p\|_V.\]This
completes the proof of the first part.

For the second part, it is enough to show that $\Psi_0(z)= W =0$.
For this, we show that $A^*$ is a completely non-unitary operator.
To this end, first notice that $\clq\subseteq H^2_{H_1}(\D)$.
Indeed, for each $f\in H^2_{H_0}(\D)$ we have\[f_n := M_{\Psi}^{*n}
f = W^{*n} f \in H^2_{H_0}(\D) \quad \quad (n  \in \Nat),\]and hence
for $g \in \clq$ we have
\begin{align*}
|\langle f, g\rangle| = |\langle M_{\Psi}^n f_n,g\rangle| = |\langle
f_n, M_{\Psi}^{* n}g\rangle| = |\langle f_n, T_2^{* n}g\rangle| \leq
\|f_n\| \|T_2^{*n} g\| = \|f\| \|T_2^{*n} g\|  \quad (n\in\Nat).
\end{align*}
Since $T_2$ is a pure contraction, $\langle f, g\rangle=0$. This
implies that $\clq\subseteq H^2_{H_1}(\D)$.

\NI On the other hand, note that $M_z $
is the minimal isometric dilation of $T_1$, that is,
\[\bigvee_{n\ge 0} M_z^n \clq=
H^2_{\cld_{T_1}}(\D),\]and $H^2_{H_1}(\D)$ is a $M_z$-reducing subspace of $H^2_{\cld_{T_1}}(\D)$. Therefore
\[H^2_{H_0}(\D)=\{0\}.\]
This shows $H_0 = \{0\}$ and completes the proof of this theorem.
\qed

In the special case where $\dim \clh < \infty$ and where $(T_1,
T_2)$ is a commuting pair of pure contractions, we have
\[\sigma(T_j) \cap
\partial \D = \emptyset \quad \quad (j = 1, 2).\]
Hence, in this particular case, we recover Agler and McCarthy's
sharper von Neumann inequality for commuting pairs of strictly
contractive matrices (see Theorem 3.1 in \cite{AM1}). Moreover, the
present proof is more direct and explicit than the one by Agler and
McCarthy (see, for instance, case (ii) in page 145 \cite{AM1}).

\newsection{Concluding remarks}\label{S-c}

\NI\textsf{Uniqueness of varieties:} Let $(T_1, T_2)$ be a pair of
pure commuting contractions on a Hilbert space $\clh$ and $\dim
\cld_{T_j} < \infty$, $j = 1, 2$. Theorem \ref{pure isometric
dilation} implies that there exists a rational inner function
${\Psi} \in H^\infty_{\clb(\cld_{T_1})}(\D)$ such that
\[{\Pi} T_1^* = M_z^* {\Pi}, \quad \mbox{and} \quad {\Pi} T_2^*
= M_{{\Psi}}^* {\Pi},\]where $\Pi : \clh \raro H^2_{\cld_{T_1}}(\D)$
is the minimal isometric dilation of $T_1$ (see (\ref{dil-def}) in
Section \ref{prel}). Furthermore, by the second assertion of Theorem
\ref{vn-new1}, we have \[\|p(T_1, T_2)\|_{\clh} \leq \|p\|_{{V}}
\quad \quad (p \in \mathbb{C}[z_1, z_2]),\]where
\[{V} = \{(z_1, z_2) \in \D^2: \det({\Psi}(z_1)-z_2
I)=0\}\]is a distinguished variety.

Now let $\tilde{\Pi} : \clh \raro H^2_{\cld_{T_2}}(\D)$ be the
minimal isometric dilation of $T_2$. Then Theorem \ref{pure
isometric dilation} applied to $(T_2, T_1)$ yields a rational inner
multiplier $\tilde{\Psi} \in H^\infty_{\clb(\cld_{T_2})}(\D)$ such
that
\[\tilde{\Pi} T_2^* = M_z^* \tilde{\Pi} , \quad \mbox{and}
\quad\tilde{\Pi} T_1^* = M_{\tilde{\Psi}}^* \tilde{\Pi}.\]Therefore
$(M_{\tilde{\Psi}}, M_z)$ on $H^2_{\cld_{T_2}}(\D)$ is also an
isometric dilation of $(T_1, T_2)$. Furthermore, \[\Psi =
\tau_{U^*}, \quad \mbox{and} \quad \tilde{\Psi} =
\tau_{\tilde{U}},\]where \[U = \begin{bmatrix} A & B
\\C& D\end{bmatrix} , \quad \mbox{and} \quad \tilde{U} = \begin{bmatrix} D^* & B^*\\ C^*&
A^*\end{bmatrix}.\] Again applying the second assertion of Theorem
\ref{vn-new1} to $(T_2, T_1)$, we see that \[\|p(T_1, T_2)\|_{\clh}
\leq \|p\|_{\tilde{V}} \quad \quad (p \in \mathbb{C}[z_1,
z_2]),\]where \[\tilde{V} = \{(z_1, z_2) \in \D^2:
\det(\tilde{\Psi}(z_2)-z_1 I)=0\}\]is a distinguished variety. Now
it follows from Lemma 1.7 in \cite{AM1} that \[V = \tilde{V}.\]

\vspace{0.2in}

\NI\textsf{Joint eignespaces:} Let $(T_1, T_2)$ and the variety $V$
be as in Theorem \ref{vn-new1}. Then the joint eigenspace of $(T_1^*,
T_2^*)$ is contained in the distinguished variety $V$. Indeed, let
\[T_j^* v=\bar{\lambda}_j v \quad \quad (j = 1, 2),\]
for some $(\lambda_1,\lambda_2)\in \D^2$ and for some non-zero
vector $v \in \clh$. Then $(\lambda_1,\lambda_2)\in V$. Equation
\eqref{U-h} gives
\[
 U(D_{T_1}v, \bar{\lambda}_1D_{T_2}v)=(\bar{\lambda}_2D_{T_1}v,D_{T_2}v).
\]
Hence, by (\ref{unitary})
\[\begin{bmatrix}A & B\\C& D \end{bmatrix} \begin{bmatrix}D_{T_1}v \\ \bar{\lambda}_1
D_{T_2}v \end{bmatrix} = \begin{bmatrix}\bar{\lambda}_2 D_{T_1}v \\
D_{T_2}v \end{bmatrix}.\]
Then by Lemma~\ref{identity},
\[(A + \bar{\lambda}_1 B (I- \bar{\lambda}_1 D)^{-1}C) (D_{T_1}v) = \bar{\lambda}_2 (D_{T_1}v),\]
and hence
\[
\left(A^*+\lambda_1C^*(I-\lambda_1D^*)^{-1}B^*\right)^*D_{T_1}v =\bar{\lambda_2} D_{T_1}v,
\]
that is, $\textup{det}(\Psi(\lambda_1)-\lambda_2I)=0$. Thus the claim
follows.

\vspace{0.2in}

\NI\textsf{An example:} We conclude this paper by pointing out a
simple but verification example of our sharper von Neumann
inequality. Let $M_z$ be the shift operator on
$H^2_{\mathbb{C}^m}(\D)$ and
\[(T_1, T_2) = (M_z,M_z).\]
Since $D_{M_z} = P_{\mathbb{C}^m}$, by a simple calculation it
follows that
the unitary $U$ in (\ref{unitary}) has the form \[U = \begin{bmatrix} 0 & W \\
I_{\mathbb{C}^m} &0
\end{bmatrix},\]
where $W\in \mathcal{B}(\mathbb{C}^n)$ is an arbitrary unitary operator.
Let us choose a unitary $W$ in $\mathcal{B}(\mathbb{C}^m)$.
In this case,
\[\tau_{U^*}(z) = \Psi(z) = z W^* \quad\quad (z \in \D).\]
Let $\{\lambda_1,\dots,\lambda_k\}$, $ 1 \leq k \leq m$, be the set
of distinct eigenvalues of $W^*$ and
\[p(z_1,z_2):= \prod_{i=1}^k(z_2-\lambda_iz_1).\] Then the
distinguished variety $V$ in the second assertion of Theorem
\ref{vn-new1} is given by \[\begin{split}V & = \{(z_1, z_2)\in\D^2 :
\det (z_1 W^* - z_2 I_{\mathbb{C}^m}) =0\} \\ & = \{(z_1,
z_2)\in\D^2 : p(z_1,z_2)=0\},\end{split}\]and hence for any $p \in
\mathbb{C}[z_1, z_2]$, the inequality
\[\|p(M_z, M_z)\|_{\clb(H^2_{\mathbb{C}^m}(\D))} \leq \|p\|_V\]holds.
In particular, if we choose $W=I_{\mathbb{C}^m}$, then
the distinguished variety $V$ is given by
\[V = \{(z,z): z\in\D\}.\] This observation also follows by a direct
calculation.

\end{document}